\documentclass{amsart}

\usepackage{amsfonts}
\usepackage{amsmath}
\usepackage{amsthm}
\usepackage{amssymb}
\usepackage[english]{babel}
\usepackage{graphics}
\usepackage{latexsym}
\usepackage{longtable} \usepackage{mathrsfs}
\usepackage{graphicx}
\usepackage{wrapfig} 
\usepackage[pdftex,colorlinks=true]{hyperref}

\newtheorem{theorem}{Theorem}
\newtheorem{corollary}[theorem]{Corollary}
\newtheorem{lemma}[theorem]{Lemma}
\newtheorem{proposition}[theorem]{Proposition}

\newtheorem{claim}[theorem]{Claim}

\theoremstyle{definition}
\newtheorem{definition}[theorem]{Definition}

\newcommand{\mD}{\mathcal{D}}

\newcommand{\R}{\mathbb{R}}
\newcommand{\N}{\mathbb{N}}

\newcommand{\Z}{\textbf{Z}}
\renewcommand{\P}{\textrm{P}}

\newcommand{\noi}{\noindent}
\newcommand{\ms}{\medskip}

\newcommand{\al}{\alpha}
\newcommand{\be}{\beta}

\newcommand{\De}{\Delta}

\newcommand{\Om}{\Omega}

\newcommand{\larrow}{\longrightarrow}
\newcommand{\ot}{\otimes}

\newcommand{\ri}{\rightarrow}

\newcommand{\sub}{\subseteq}

\newcommand{\ess}{\textrm{ess}}

\newcommand{\Div}{\textrm{Div}}

\newcommand{\bt}{\begin{theorem}}\newcommand{\et}{\end{theorem}}
\newcommand{\bd}{\begin{definition}}\newcommand{\ed}{\end{definition}}
\newcommand{\bl}{\begin{lemma}}\newcommand{\el}{\end{lemma}}
\newcommand{\beq}{\begin{equation}}\newcommand{\eeqq}{\end{equation}}
\newcommand{\bc}{\begin{claim}}\newcommand{\ec}{\end{claim}}
\newcommand{\bp}{\begin{proof}}\newcommand{\ep}{\end{proof}}

\newcommand{\BPP}{\medskip \noindent \textbf{Proof of Proposition} }
\newcommand{\BPT}{\medskip \noindent \textbf{Proof of Theorem} }

\numberwithin{equation}{section}
\numberwithin{theorem}{section}

\begin{document}

\title[A H\"older Continuous Nowhere Improvable Function with Derivative...]{A H\"older Continuous Nowhere Improvable Function with Derivative Singular\\ Distribution}
\author{\textsl{Nikolaos I. Katzourakis}}
\address{BCAM - Basque Center for Applied Mathematics, Alameda de Mazarredo 14, E-48009, Bilbao, Spain}
\email{nkatzourakis@bcamath.org}


\date{}



\maketitle

\begin{abstract} We present a class of functions $\mathcal{K}$ in $C^0(\R)$ which is variant of
the Knopp class of nowhere differentiable functions. We
derive estimates which establish $\mathcal{K} \sub C^{0,\al}(\R)$ for $0<\al<1$ but no $K \in \mathcal{K}$ is pointwise anywhere improvable to $C^{0,\be}$ for any $\be>\al$. In particular, all $K$'s are nowhere differentiable with derivatives singular distributions. $\mathcal{K}$ furnishes explicit realizations of the functional analytic result of Berezhnoi \cite{Be}.

Recently, the author and simulteously others laid the foundations of Vector-Valued Calculus of Variations in $L^\infty$ \cite{K2, K3, K4}, of $L^\infty$-Extremal Quasiconformal maps \cite{CR, K5} and of Optimal Lipschitz Extensions of maps \cite{SS}. The  ``Euler-Lagrange PDE'' of Calculus of Variations in $L^\infty$ is the nonlinear nondivergence form Aronsson PDE with as special case the $\infty$-Laplacian. 

Using $\mathcal{K}$, we construct singular solutions for these PDEs. In the scalar case, we partially answered the open $C^1$ regularity problem of Viscosity Solutions to Aronsson's PDE \cite{K1}. In the vector case, the solutions can not be rigorously interpreted by existing PDE theories and justify our new theory of Contact solutions for fully nonlinear systems \cite{K6}. Validity of arguments of our new theory and failure of classical approaches both rely on the properties of $\mathcal{K}$.
\end{abstract}

\section{Introduction}

\noi Let $\al \in (0,1)$ and $\nu \in \N$ be fixed parameters. We
define the continuous function $K_{\al,\nu} : \R \larrow [0,1]$ by
 \beq \label{eq1}
 K_{\al,\nu}(x)\ :=\ \sum_{k=0}^\infty 2^{-2\al \nu k}\phi(2^{2 \nu
 k}x),
 \end{equation}
where $\phi$ is a sawtooth function, given by $\phi(x):=|x|$ when $x
\in [-1,1]$ and extended on $\R$ as a  periodic function by setting
$\phi(x+2):=\phi(x)$. Explicitly,
 \begin{equation} \label{eq2}
 \phi(x)\ =\ \sum_{i=-\infty}^{+\infty}\big|x-2i\big|\chi_{(i-1,i+1]}(x).
 \end{equation}
Formulas \eqref{eq1}, \eqref{eq2} introduce a parametric family in
the space of H\"older continuous functions $C^{0,\al}(\R)$
which are not differentiable at any point of $\R$. The first
examples of nowhere differentiable continuous functions given by
Weierstrass, Bolzano and Cell\'erier have been followed by numerous
functions well behaved with respect to continuity but very singular
with respect to differentiability. Our example of $K_{\al,\nu}$ is a
variant of the Knopp function \cite{Kn} (see also \cite{B-D} and
\cite{C}) and relates directly to several other examples existing in
the literature, for example the Takagi-Van der Waerden functions, as
well as the McCarthy function \cite{M}.

Our explicit class of functions gives a simple realization of the abstract functional analytic result of Berezhnoi \cite{Be}, who proved that every infinite-dimensional Banach space of functions which enjoy some degree of regularity, contains an infinite-dimensional closed subspace of functions ``nowhere improvable", namely not smoother than the least smooth function in the space.

It is worth-noting that examples of continuous nowhere differentiable functions
still attract mathematical interest. Recently, Allart and Kawamura
\cite{Al-K} characterized the sets at which ``improper'' infinite
derivatives exist for the Takagi function, while Lewis \cite{L}
studies probabilistic aspects of the Katsuura function. For a
thorough historical review and an extended list of references, we
refer to Thim \cite{T}.

Herein we derive suitable precise estimates which establish that
$K_{\al,\nu}$ is in the H\"older space $C^{0,\al}(\R)$ for all $\nu
\in \N$, but if $\nu$ is sufficiently large ($2\nu>1/(1-\al)$) the function is at no point improvable to a H\"older continuous $C^{0,\be}$ function for any $\be \in(\al,1]$. \emph{A function $f\in C^0(\R)$ is called H\"older continuous $C^{0,\be}$ at $x \in \R$ if there exist $r,C>0$ such that}
\begin{equation}
|f(y)-f(x)|\ \leq\ C|y-x|^\be,
\end{equation}
\emph{for all $y\in[x-r,x+r]$}. As a consequence, for $\be=1$ we deduce that $K_{\al,\nu}$ is nowhere differentiable since the pointwise derivative does not exist anywhere.

Since BV functions are differentiable almost everywhere, the distributional derivative $D K_{\al,\nu}$ in
$\mD'(\R)$ is a singular first order distribution and can
not be realized by a signed measure.  In particular, for any $\al,\be \in (0,1)$, $\be>\al$ and $\nu$ large enough, the $C^{0,\al}$-function $K_{\al,\nu}$ is neither a $BV$ nor a $C^{0,\be}$
function no matter how ``close'' to the Lipschitz space
$C^{0,1}(\R)$ it might be.
\[
\label{Fig. 1} \underset{\text{Fig.\ 1: Simulation of $K_{\al,\nu}$
over $(0,2)$ with $\nu=2$, $\al=1/2$ (50 terms,
Mathematica)}}{\centerline{\includegraphics[scale=0.15]{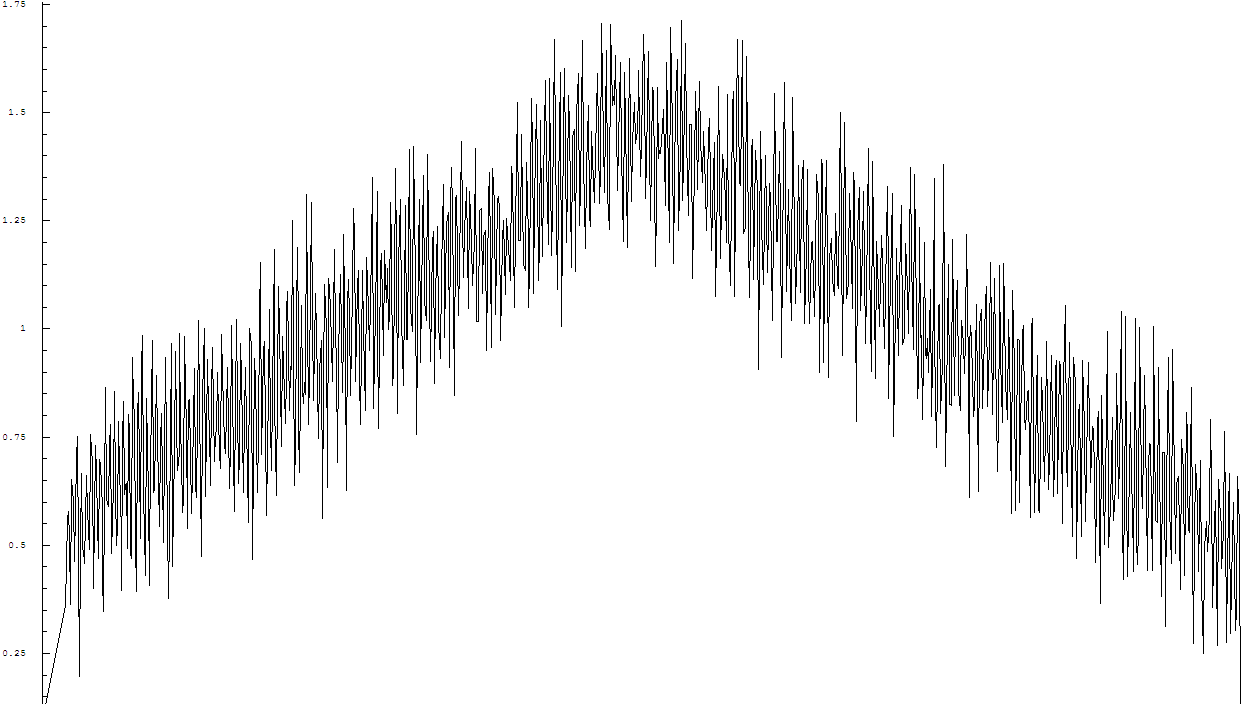}}}
\]
\[
\underset{\text{Fig.\ 2: Simulation of $K_{\al,\nu}$
over $(0,2)$ with $\nu=4$, $\al=1/8$ (50 terms,
Mathematica)}}{\centerline{\includegraphics[scale=0.15]{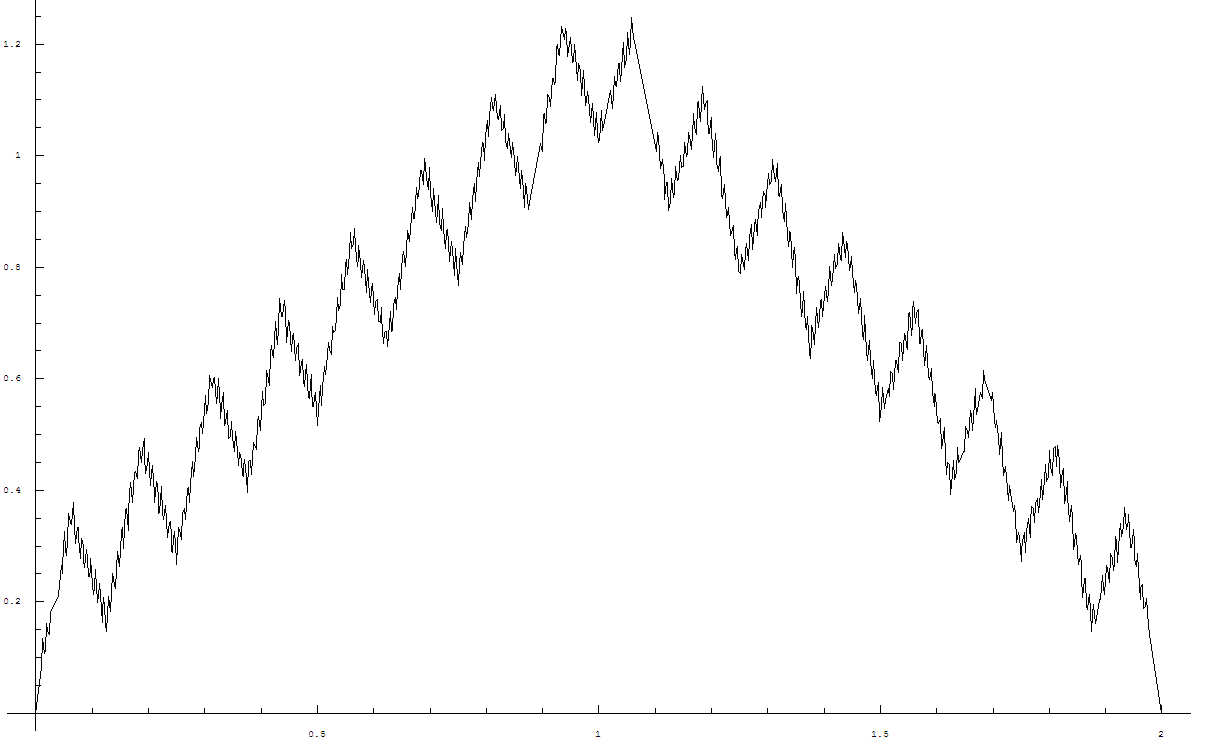}}}
\]
\[
\label{Fig. 3} \underset{\text{Fig.\ 3: Simulation of $K_{\al,\nu}$
over $(0,2)$ with $\nu=4$, $\al=5/8$ (50 terms,
Mathematica)}}{\centerline{\includegraphics[scale=0.15]{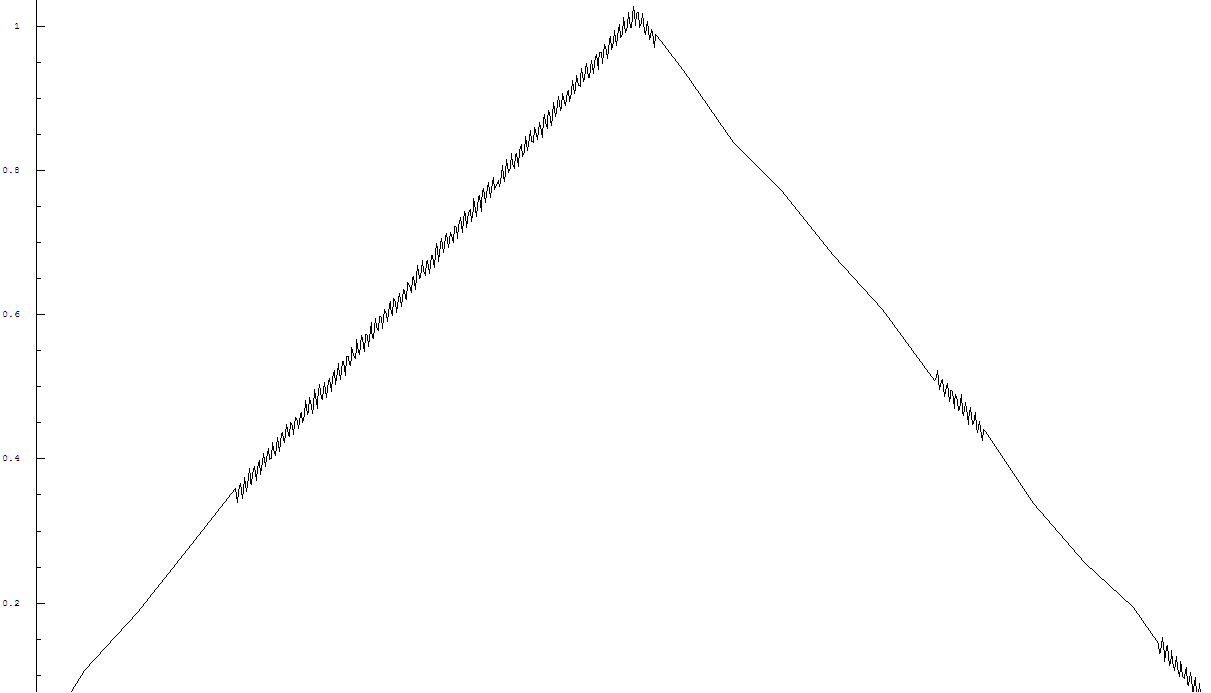}}}
\]
The necessity to construct pathological functions which bear the
specific properties of $K_{\al,\nu}$ originates from the theory of
nonlinear partial differential equations, especially regularity
theory of degenerate 2nd order elliptic partial differential
equations and systems including the celebrated $\infty$-Laplacian
 \begin{equation}
 \label{eq3a}
 \Delta_{\infty}u\ :=\  Du \otimes Du :D^2 u\ = \ 0, \ \ u \ :
 \R^n \larrow \R
 \end{equation}
(that is $\Delta_{\infty}u= D_i u \, D_j u\, D^2_{ij}u$, with the summation convention employed),
as well as the more general Aronsson equation
 \begin{equation}
 \label{eq3b}
 A_\infty u\ :=\ H_p(Du) \otimes H_p(Du) :D^2 u\ = \ 0,
 \end{equation}
for a Hamiltonian $H \in C^1(\R^n)$ and $H_p(p):= DH(p)$. The Aronsson PDE is the ``Euler-Lagrange PDE'' of Calculus of Variations in the space $L^\infty$ for the supremal functional
\begin{equation}
E_\infty(u,\Om)\ :=\ \underset{\Om} {\ess \, \sup}\, H(Du), \ \ \ u \in W^{1,\infty}(\Om).
\end{equation}
The celebrated $\infty$-Laplacian corresponds to the model functional $\|Du\|_{L^\infty(\Om)}$ when we chose as Hamiltonian $H$ the Euclidean norm. The name ``$\infty$-Laplacian'' originates in its first derivation in the limit of the $p$-Laplacian 
\begin{equation}
\De_p u\ :=\ \Div(|Du|^{p-2}Du)\ =\ 0
\end{equation}
as $p\ri \infty$ by Aronsson. The $p$-Laplacian is the Euler-Lagrange PDE of the $p$-Dirichlet functional
\begin{equation}
E_p(u,\Om)\ :=\ \int_{\Om} |Du|^p, \ \ \ u \in W^{1,p}(\Om).
\end{equation}
When passing to the limit $p\ri \infty$, divergence structure is lost and, unlike $\De_p$, the operator $\De_\infty$ is quasilinear but not in nondivergence form. Hence, standard weak and distributional solution approaches of modern PDE theory do not work. In \cite{Ar1} and \cite{Ar2} Aronsson constructed singular solutions to $\Delta_{\infty}u=0$, while the general $C^1$ regularity problem related to $\Delta_{\infty}$ and $A_\infty$ is still open, except for the dimension
$n=2$ (\cite{S, E-S, J-N, W-Y, Cr2}). In the Author's work \cite{K1}, by employing the function of this paper we gave a partial negative answer to this conjecture by showing that there exist Hamiltonians for which the Aronsson PDE admits non-$C^1$ solutions. 

The general vector case of the $\infty$-Laplacian for maps is much more intricate and its study started only recently in \cite{K2, K3, K4, K5}, where the foundations of Vector-Valued Calculus of Variations in $L^\infty$ and its ``Euler-Lagrange PDE system'' have been laid. Related simultaneous results appeared also in \cite{CR, SS}. Capogna and Raich in \cite{CR} simultaneously but independently used as Hamiltonian the so-called trace distortion $|Du|^n/\det(Du)$ defined on local diffeomorphisms and developed a parallel to the Author's approach for Extremal $L^\infty$-Quasiconformal maps. Also, Sheffield and Smart in \cite{SS}, developed the related subject of Vector-Valued Lipschitz Extensions by using as Hamiltonian the operator norm $\| . \|$ on the space $\R^N \ot \R^n$ of gradients $Du$ of maps $u : \R^n \larrow \R^N$.  

For general maps $u:\R^n \larrow \R^N$, the $\infty$-Laplacian in the vector case reads
\begin{equation} \label{inftylap}
\De_\infty u\ :=\ Du \ot Du :D^2u \ +\ |Du|^2[Du]^\bot\De u\ = \ 0.
\end{equation}
Here $[Du(x)]^\bot$ is the projection on the nullspace $N(Du(x)^\top)$ of the (transpose of the gradient) operator $Du(x)^\top : \R^N \larrow \R^n$. In index form reads
\begin{equation} 
D_i u_\al \, D_j u_\be\, D_{ij}^2u_\be \ +\ |Du|^2[Du]_{\al \be}^\bot D^2_{ii} u_\be\ = \ 0
\end{equation}
and was first derived in \cite{K2}. The general Aronsson PDE system corresponding to a rank-one convex Hamiltonian $H\in C^2(\R^N \ot \R^n)$ is
\begin{equation} \label{aronsyst}
A_\infty u\ :=\ \Big(H_P\ot H_P\, +\,  H[H_P]^\bot H_{PP}\Big)(Du):D^2u\ = \ 0.
\end{equation}
Here $[H_P(Du(x))]^\bot$ is the projection on the nullspace $N(H_P(Du(x))^\top)$ of the operator $H_P(Du(x))^\top : \R^N \larrow \R^n$. For details we refer to Section \ref{PDE} and \cite{K2}. 

The vector case of \eqref{inftylap} and \eqref{aronsyst} is extremely difficult. The main reason is the existence of singular ``solutions" constructed by means of the functions in this paper which show that under the current state-of-art in PDE theory, such systems can not be studied rigorously and we can not even interpret appropriately their singular solutions. It is a similar problem to that of interpretation of the Dirac $\delta$ in Quantum Mechanics before measure theory. 

Moreover, \eqref{inftylap} and \eqref{aronsyst} are nonlinear, nonmonotone and in nondivergence form and have discontinuous coefficients even for $C^\infty$ solutions: the normal projections $[Du]^\bot$ is discontinuous when the rank of $Du$ changes. This is a genuinely vectorial phenomenon and happens because there exist smooth $\infty$-Harmonic maps whose rank of the gradient is not constant: such an example is given by 
\begin{equation} \label{1.13}
u(x,y)\, :=\, e^{ix}-e^{iy} \ ,\ \ \ u \ :\ \{|x \pm y|<\pi\}\sub \R^2 \larrow \R^2. 
\end{equation}
Indeed, \eqref{1.13} is $\infty$-Harmonic on the rhombus and has $rk(Du)=1$ on the diagonal $\{x=y\}$, but has $rk(Du)=2$ otherwise and the projection $[Du]^\bot$ is discontinuous (for more details see \cite{K2}). 
\[
\begin{array}{c}
\includegraphics[scale=0.15]{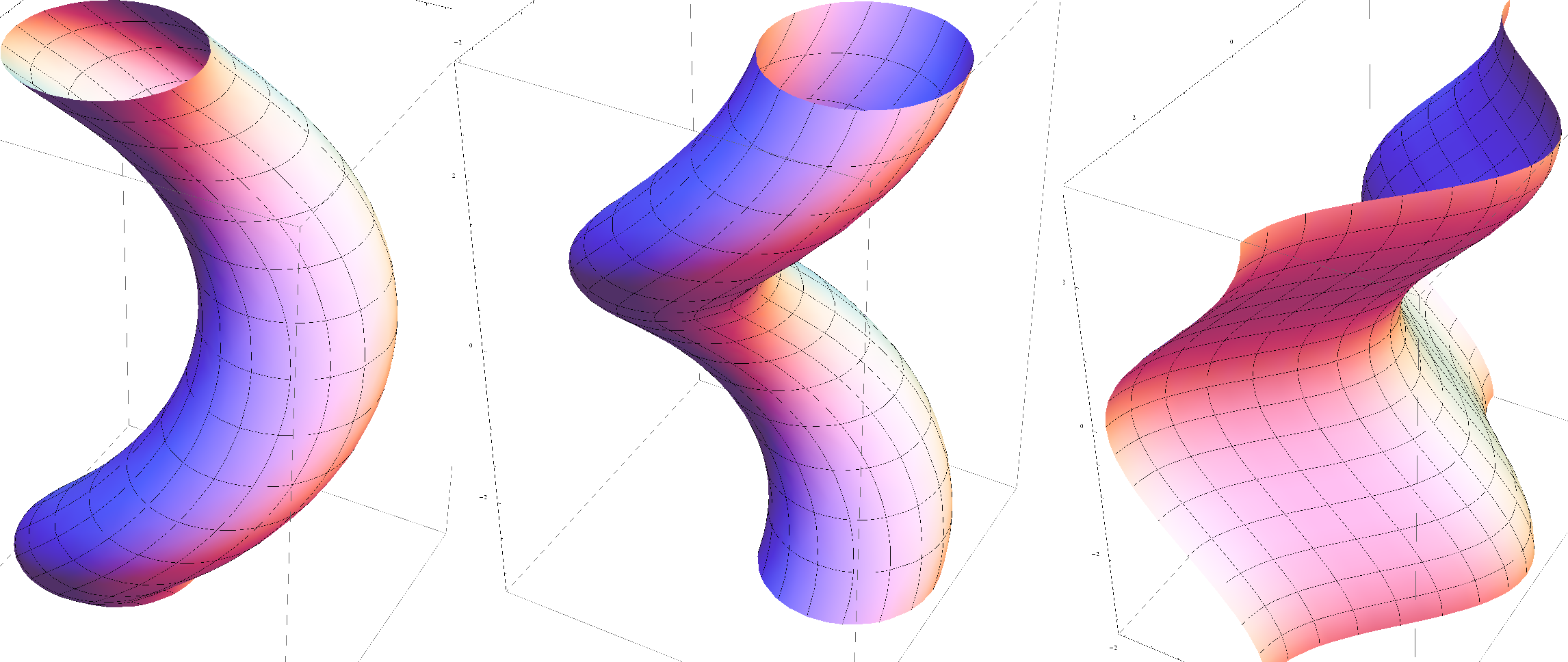}\\
\underset{\text{3 projections on $\R^3$ of the graph of $u(x,y)=e^{ix}-e^{iy}$, its range on $\R^2$ and its ``covering sheets''.}}
{\includegraphics[scale=0.18]{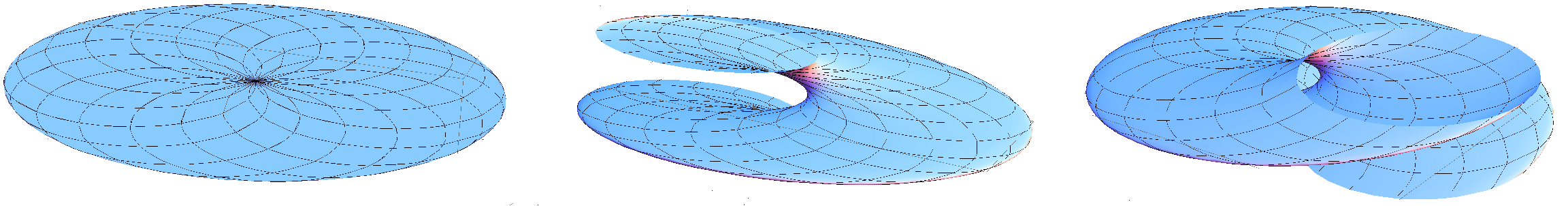}}
\end{array}
\]
In general, \emph{$\infty$-Harmonic maps present a phase separation},  with a certain hierarchy. On each phase the dimension of the tangent space is constant and these phases are separated by \emph{interfaces} whereon the rank of $Du$ ``jumps'' and $[Du]^\bot$ gets discontinuous. 

The related problems of the scalar case were unsolved for some years and were finally settled in the '90s with the advent of Viscosity Solutions. However, viscosity solutions apply only to scalar PDEs and monotone PDE systems. In \cite{K6} we have intoduced an appropriate systematic theory which applies to fully nonlinear PDE systems and in particular allows to study \eqref{inftylap} and \eqref{aronsyst} rigorously and effectively. This theory extends Viscosity Solutions in the vector case of systems and is based on the discovery of an extremality principle which applies to maps. Contact Solutions bear stability properties similar to their scalar counterparts of Viscosity Solutions and this feature renders them extremely efficient when trying to prove existence results via approximation.

This paper is organized as follows: in Section \ref{K} we present the basic material about our singular class of functions and in Section \ref{PDE} we present some material related to singular PDE solutions.

\section{The Singular Function $K$.} \label{K}

The following Theorem lists the properties of $K_{\al,\nu}$.

\begin{theorem} \label{th1} (i) The function $K_{\al,\nu}$ is in $C^{0,\al}(\R)$ for all
$\nu \in \N$. Moreover, we have the uniform bound $0\leq K_{\al,\nu}
\leq 1/(1-2^{-2\nu\al})$ and if  $x,y \in \R$ with $|x-y|\leq2$, we have the estimate
 \begin{equation} \label{eq4}
 \big|K_{\al,\nu}(x)-K_{\al,\nu}(y) \big|\ \leq \
  C(\al,\nu)\ |x-y|^\al,
 \end{equation}
where
 \begin{equation} \label{eq5}
C(\al,\nu)\ :=\ \frac{1}{1-2^{-2\nu(1-\al)}}
  + \frac{2}{2^{2\nu(\al-1)}-2^{-2\nu}}.
  \end{equation}

\noi (ii) If $\al\in (0,1)$ and $2\nu>1/(1-\al)$, then
$K_{\al,\nu}$ is pointwisely nowhere improvable to $C^{0,\be}$ on $\R$ for any  $\be \in (\al,1]$. Moreover, for any $x\in \R$, $m \in \N$, we have
 \begin{equation} \label{eq6}
\frac{\big|K_{\al,\nu}(x+t_m(x))-K_{\al,\nu}(x)\big|}{|t_m(x)|^\be}\
\geq\ \ K(m,\nu,\al,\be),
  \end{equation}
where
 \begin{equation} \label{eq7}
 K(m,\nu,\al,\be) \ :=\ \frac{2^{\be-1}(2^{2\nu(1-\al)} - 2)}
{2^{2\nu(1-\al)} - 1} \big(2^{2\nu(\be-\al)}\big)^m
 \end{equation}
and $t_m :\R \larrow \{\pm 2^{-2\nu m -1}\}$ is the step function
given by
 \begin{equation} \label{eq8}
t_m (x) \ := \ 2^{-2\nu m -1} \sum_{i=-\infty}^{+\infty}\Big[
\chi_{(i,i+ \frac{1}{2}]}\big(2^{2\nu m}x\big)\ -\
\chi_{(i+\frac{1}{2},i+1]}\big(2^{2\nu m}x\big) \Big].
 \end{equation}.
\end{theorem}

As noted earlier, nowhere differentiable functions have genuine distributional derivatives, since, BV functions must necessarily be differentiable almost everywhere. As a corollary of the previous theorem, we provide a lower bound on the total variation of the difference quotients which establishes this fact without employing the fine properties of BV functions.

\begin{proposition} \label{pr1} If $\al\in (0,1)$ and $2\nu>1/(1-\al)$, for any $M\geq 1$,
$m \in \N$, we have the following lower bounds in $L^1_{loc}(\R)$ for the difference quotients  
 \begin{equation} \label{eq10}
\frac{1}{2M}\int_{-M}^{+M}\left|\frac{K_{\al,\nu}(x+2^{-2\nu m
-1})-K_{\al,\nu}(x)}{2^{-2\nu m -1}}\right|dx\ \geq\ \
\frac{1}{4}K(m,\nu,\al,1).
  \end{equation}
\end{proposition}

\ms 

The previously obtained estimates readily imply the following

\ms 

\begin{corollary} \label{cor2}
For any $\al \in (0,1)$ and $2\nu>1/(1-\al)$ with $\nu \in \N$, $x\in \R$, $\be \in (\al,1]$ and $M\geq1$, we have
 \begin{equation} \label{eq9}
\limsup_{t\rightarrow
0}\frac{\big|K_{\al,\nu}(x+t)-K_{\al,\nu}(x)\big|}{|t|^\be}\ = \
+\infty.
 \end{equation}
Hence, the $C^{0,\al}$ function $K_{\al,\nu}$ is nowhere improvable to $C^{0,\be}$. In particular, if $\be=1$ then the function is pointwisely nowhere differentiable on $\R$. Also,
 \begin{equation} \label{eq11}
\limsup_{t\rightarrow 0} \frac{1}{2M}\int_{-M}^{+M}
\left|\frac{K_{\al,\nu}(x+t)-K_{\al,\nu}(x)}{t}\right|dx\ = \
+\infty.
 \end{equation}
Hence, the difference quotients are unbounded in $L^1_{loc}(\R)$ and the distributional derivative of $K_{\al,\nu}$ does not exist as a signed measure.
\end{corollary}

The first part of Corollary \ref{cor2} is immediate, while the second follows by estimate \eqref{eq10} and application of the folklore fact that an $L^1_{loc}(\R)$ function is of Bounded Variation if and only if the difference quotients converge weakly$^*$ in the sense of measures.

\BPT \ref{th1}. (i) We begin by observing that \eqref{eq2} implies
$|\phi|\leq 1$ and hence the bound
 \begin{equation} \label{eq12}
 0\ \leq \ K_{\al,\nu} \ \leq \ \sum_{k=0}^\infty \big(2^{-2\al \nu}\big)^k \ =\
 \frac{1}{1-2^{-2\nu\al}}.
 \end{equation}
Let now $p,q \in \N$ with $p<q$ and $x \in \R$. Again by
\eqref{eq2}, we have
 \begin{equation} \label{eq13}
 \left| \sum_{k=0}^q 2^{-2\al \nu k}\phi(2^{2 \nu
 k}x) - \sum_{k=0}^p 2^{-2\al \nu k}\phi(2^{2 \nu
 k}x) \right| \ \leq \ \sum_{k=p+1}^q \big(2^{-2\al \nu}\big)^k,
 \end{equation}
which tends to $0$ as $p,q \larrow \infty$. By \eqref{eq13},
\eqref{eq1} defines a continuous function: $K_{\al,\nu} \in
C^{0}(\R)$. Fix now $x,y$ in $\R$ with $x\neq y$ and choose $t\geq1$
and $p \in \N$ such that
 \begin{equation} \label{eq14}
 |x|,\ |y| \ \leq \ 2^{2\nu-1}t
 \end{equation}
and
 \begin{equation} \label{eq15}
\frac{t}{2^{2\nu p}} \ \leq \ |y-x| \ \leq \
 \frac{t}{2^{2\nu (p-1)}}.
 \end{equation}
Since by \eqref{eq2} $\phi$ is non-expansive, that is
$|\phi(t)-\phi(s)|\leq |t-s|$ and also $|\phi|\leq 1$, we estimate
 \begin{align} \label{eq16}
\frac{\big|K_{\al,\nu}(x)-K_{\al,\nu}(y) \big|}{|x-y|^\al} \ \leq &
\ |x-y|^{-\al} \left[\sum_{k=0}^{p-1}2^{-2\al\nu k}\big| \phi(2^{\nu
k}x)-\phi(2^{\nu k}y)\big| \ + \ 2\sum_{k=p}^\infty 2^{-2\al \nu
k}\right]\nonumber\\
\leq & \ |x-y|^{-\al} \left[\sum_{k=0}^{p-1}2^{2\nu k(1-\al)}|x-y| \
+ \ 2\sum_{k=p}^\infty 2^{-2\al \nu k}\right].
\end{align}
Hence, by \eqref{eq15}, estimate \eqref{eq16} gives
\begin{align} \label{eq17}
\frac{\big|K_{\al,\nu}(x)-K_{\al,\nu}(y) \big|}{|x-y|^\al} \ \leq &
\ |x-y|^{-\al} \left[\sum_{k=0}^{p-1}2^{2\nu k (1-\al)}|x-y| \ +
2\frac{2^{-2\nu \al p}}{1-2^{-2\nu \al}}\right] \nonumber\\
 = & \ |x-y|^{-\al}
\left[\frac{2^{2\nu p(1-\al)} -1}{2^{2\nu (1-\al)} -1}|x-y| \ +
2\frac{2^{-2\nu \al p}}{1-2^{-2\nu \al}}\right] \\
 \leq & \ |x-y|^{-\al}
\left[\frac{2^{2\nu p(1-\al)} -1}{2^{2\nu (1-\al)} -1}|x-y| \ +
2\frac{2^{-2\nu \al p}}{1-2^{-2\nu \al}}
\frac{2^{2\nu p}}{t}|x-y|\right]  \nonumber\\
 \leq & \ |x-y|^{1-\al}
\left[\frac{2^{2\nu p(1-\al)} -1}{2^{2\nu (1-\al)} -1} \ +
2\frac{2^{2\nu p (1-\al)}}{1-2^{-2\nu \al}}\right]. \nonumber
\end{align}
Again by \eqref{eq15}, estimate \eqref{eq17} gives
\begin{align} \label{eq18}
\frac{\big|K_{\al,\nu}(x)-K_{\al,\nu}(y) \big|}{|x-y|^\al} \
 \leq
 &   \ \big(t2^{-2\nu (p-1)}\big)^{1-\al}2^{2\nu p(1-\al)}
\left[\frac{1-2^{-2\nu p(1-\al)}}{2^{2\nu (1-\al)} -1} \ +
\frac{2}{1-2^{-2\nu \al}}\right] \nonumber\\
=  &   \ t^{1-\al}  2^{2\nu (1-\al)} \left[\frac{1-2^{-2\nu
p(1-\al)}}{2^{2\nu (1-\al)} -1} \ +
\frac{2}{1-2^{-2\nu \al}}\right]\\
\leq  &   \ t^{1-\al}  2^{2\nu (1-\al)} \left[\frac{1}{2^{2\nu
(1-\al)} -1} \ + \frac{2}{1-2^{-2\nu \al}}\right].\nonumber
\end{align}
By estimate \eqref{eq14}, we have $t \geq
\max\{1,2^{1-2\nu}|x|,2^{1-2\nu}|y|\}$. By minimizing \eqref{eq18}
with respect to all such $t$'s, we obtain
 \begin{equation}  \label{eq19}
\frac{\big|K_{\al,\nu}(x)-K_{\al,\nu}(y) \big|}{|x-y|^\al} 
 \leq 
\left(\max\{|x|,|y|,2^{2\nu}\}\right)^{1-\al}\left[\frac{1}{2^{2\nu(1-\al)}-1}
  + \frac{2}{1-2^{-2\nu\al}}\right].
 \end{equation}
By periodicity of $K_{\al,\nu}$, we may further assume that $|x|$, $|y| \leq 1$. Hence, estimate \eqref{eq19} leads directly to \eqref{eq4} and \eqref{eq5}.

$$ $$

\noi (ii) Fix $x\in \R$ and $m \in \N$. Let $t_m : \R \larrow \R$ be
the step function given by formula \eqref{eq8}, which we reformulate
as
 \begin{equation} \label{eq20*}
t_m(x)\ = \left\{ \begin{array}{l} +2^{-2\nu m -1},
      \hspace{61pt} i2^{-2\nu m}  <\ x \ \leq  i2^{-2\nu m}+ 2^{-2\nu m -1},  \ i \in \Z.\\
                              -2^{-2\nu m -1},
     \ \ \  i2^{-2\nu m}+ 2^{-2\nu m -1}  < \ x \ \leq i2^{-2\nu m}+2^{-2\nu m},
      \ \ \ \ i \in \Z.\\
          \end{array}
 \right.
 \end{equation}
We observe that since $|t_m(x)|=\frac{1}{2}2^{-2\nu m}$ and
 \begin{equation} \label{eq20}
 \Big|2^{2\nu m}(x+t_m(x)) \ - \ 2^{2\nu m}x\Big| \ = \ \frac{1}{2},
 \end{equation}
$t_m$ is defined in such a way that no integer lies between $2^{2\nu
m}x$ and $2^{2\nu m}(x+t_m(x))$. By \eqref{eq1}, we can first estimate
from below the difference quotient
$\big|\big(K_{\al,\nu}(x+t_m(x))-K_{\al,\nu}(x)\big)/t_m(x)\big|$ for $\be=1$ as
 \begin{align}
  \label{eq21}
 \left| \frac{K_{\al,\nu}(x+t_m(x))-K_{\al,\nu}(x)}{t_m(x)} \right|
 \ \geq & \  \Bigg| \sum_{k=m+1}^\infty 2^{-2\al \nu k}
  \left(\frac{\phi(2^{2 \nu k}(x+t_m(x))) - \phi(2^{2 \nu k}x)}{t_m(x)}\right)\nonumber\\
  & \ \ + \ 2^{-2\al \nu m}
  \left(\frac{\phi(2^{2 \nu m}(x+t_m(x))) - \phi(2^{2 \nu m}x)}{t_m(x)}\right)
  \Bigg|\\
  &  - \ \sum_{k=0}^{m-1} 2^{-2\al \nu k} \left|
  \frac{\phi(2^{2 \nu k}(x+t_m(x))) - \phi(2^{2 \nu
  k}x)}{t_m(x)}\right|. \nonumber
  \end{align}
We will derive estimate \eqref{eq6} by estimating each term of
\eqref{eq21}. First, observe that the sum $\sum_{k=m+1}^\infty$ in
\eqref{eq21} vanishes, since by \eqref{eq2} $\phi$ is $2$-periodic:
indeed, for $k\geq m+1$, we have
 \begin{align} \label{eq23}
\phi(2^{2 \nu k}(x+t_m(x))) - \phi(2^{2 \nu k}x)   \ & = \ \phi
\big( 2^{2 \nu k}x \pm 2^{2 \nu (k-m)-1} \big) \ - \ \phi \big( 2^{2
\nu k}x\big) \nonumber\\
& = \ \phi \big( 2^{2 \nu k}x \pm 2\, 2^{2 (\nu (k-m)-1)} \big) \ - \
\phi \big( 2^{2 \nu k}x\big)\\
& = \ 0 , \nonumber
 \end{align}
the last equality being obvious since $2^{2 (\nu (k-m)-1)} \in \N$.
Next, the sum $\sum_{k=0}^{m-1}$ in \eqref{eq21} can be estimated as
\begin{align} \label{eq24}
\sum_{k=0}^{m-1} 2^{-2\al \nu k}
  \frac{\big|\phi(2^{2 \nu k}(x+t_m(x))) - \phi(2^{2 \nu
  k}x) \big|}{| t_m(x) |}\
   &  \leq \ \sum_{k=0}^{m-1} 2^{-2\al \nu k}
  \frac{\big|2^{2 \nu k}(x+t_m(x)) - 2^{2 \nu
  k}x \big|}{| t_m(x) |} \nonumber\\
  & = \ \sum_{k=0}^{m-1} 2^{2 \nu (1-\al)k}  \\
& = \ \frac{1-2^{2\nu(1-\al)m}}{1-2^{2\nu (1-\al)}}. \nonumber
 \end{align}
Finally, by the definition of $t_m$ and the fact that $\phi$ is
piecewise affine with unit slope along the segments between
integers, the remaining middle term of \eqref{eq21} gives
\begin{align} \label{eq25}
\left| 2^{-2\al \nu m}\left(
  \frac{\phi(2^{2 \nu m}(x+t_m(x))) - \phi(2^{2 \nu
  m}x)}{ t_m(x) }\right) \right|\
   &  =  2^{-2\al \nu m}
  \frac{\big|2^{2 \nu m}(x+t_m(x)) - 2^{2 \nu
  m}x|}{ |t_m(x)\big| } \\
   &  =  2^{2 \nu (1-\al)m}.\nonumber
 \end{align}
By utilizing equations \eqref{eq23}, \eqref{eq24} and \eqref{eq25},
estimate \eqref{eq21} implies
 \begin{align} \label{eq26}
\left|\frac{K_{\al,\nu}(x+t_m(x))-K_{\al,\nu}(x)}{t_m(x)}\right|\ &
\geq\ 2^{2 \nu (1-\al)m} \ -\ \frac{1-2^{2\nu(1-\al)m}}{1-2^{2\nu
(1-\al)}} \\
& \geq \ \frac{(2^{2\nu(1-\al)} - 2)}{2^{2\nu(1-\al)} - 1} \big(2^{2\nu(1-\al)}\big)^m. \nonumber
 \end{align}
By \eqref{eq26}, if $\be \in (\al,1]$ then by employing that $|t_m(x)|=2^{-2\nu m-1}$, we have
\begin{align} \label{eq26a}
\frac{\big|K_{\al,\nu}(x+t_m(x))-K_{\al,\nu}(x)\big|}{|t_m(x)|^\be}\ &
 \geq \ 2^{-(2\nu m +1)(1-\be)}\frac{(2^{2\nu(1-\al)} - 2)}{2^{2\nu(1-\al)} - 1} \big(2^{2\nu(1-\al)}\big)^m \nonumber\\
& = \ K(m,\nu,\al,\be),
\end{align}
which is equivalent to \eqref{eq6} and \eqref{eq7}. 

\qed

\BPP \ref{pr1}.  Let $M\geq1$. We fix $m\in \N$ and set
 \begin{equation} \label{eq27}
 E_m \ :=\ \big\{x \in \R \ \big| \ t_m(x)>0 \big\}.
  \end{equation}
By \eqref{eq8}, $t_m$ is a Borel measurable function and
hence $E_m$ is a Borel set. By integrating \eqref{eq26a} on
$(-\frac{M}{2},\frac{M}{2})$ for $\be=1$, we obtain
\begin{align}
 \label{eq28}
M \frac{(2^{2\nu(1-\al)} - 2)\big(2^{2\nu(1-\al)}\big)^m}{2^{2\nu(1-\al)} - 1} \  & \leq \
\int_{-\frac{M}{2}}^{\frac{M}{2}}
\left|\frac{K_{\al,\nu}(x+t_m(x))-K_{\al,\nu}(x)}{t_m(x)}\right| dx
\nonumber\\ & = \ \int_{(-\frac{M}{2},\frac{M}{2}) \cap E_m}
\left|\frac{K_{\al,\nu}(x+t_m(x))-K_{\al,\nu}(x)}{t_m(x)}\right|
dx\\
 & \ \ + \ \int_{(-\frac{M}{2},\frac{M}{2}) \setminus E_m}
\left|\frac{K_{\al,\nu}(x+t_m(x))-K_{\al,\nu}(x)}{t_m(x)}\right| dx.
\nonumber
\end{align}
Hence, by \eqref{eq27} and \eqref{eq20*}, \eqref{eq28} gives
\begin{align} \label{eq29}
M \frac{(2^{2\nu(1-\al)} - 2)\big(2^{2\nu(1-\al)}\big)^m}
{2^{2\nu(1-\al)} - 1}   & \leq  \int_{(-\frac{M}{2},\frac{M}{2})
\cap E_m} \left|\frac{K_{\al,\nu}(x+2^{-2\nu m
-1})-K_{\al,\nu}(x)}{2^{-2\nu m -1}}\right| dx \nonumber
 \\
&  + \int_{(-\frac{M}{2},\frac{M}{2}) \setminus E_m}
\left|\frac{K_{\al,\nu}(x- 2^{-2\nu m -1})-K_{\al,\nu}(x)}{-
2^{-2\nu m -1}}\right| dx. 
\end{align}
By a change of variables in the second integral, we obtain
\begin{align} \label{eq30}
 \frac{(2^{2\nu(1-\al)} - 2)
\big(2^{2\nu(1-\al)}\big)^m}{2^{2\nu(1-\al)} - 1}  \leq  \frac{1}{M} \int_{(-\frac{M}{2},\frac{M}{2}) \cap E_m} &
\left|\frac{K_{\al,\nu}(x+2^{-2\nu m -1})-K_{\al,\nu}(x)}{2^{-2\nu m
-1}}\right| dx
\\
 +  \frac{1}{M}  \int_{((-\frac{M}{2},\frac{M}{2}) \setminus E_m) -2^{-2\nu m
-1}} & \left|\frac{K_{\al,\nu}(x+ 2^{-2\nu m
-1})-K_{\al,\nu}(x)}{2^{-2\nu m -1}}\right| dx. \nonumber
\end{align}
Hence, since $M\geq 1$ and $2^{-2\nu m -1}\leq \frac{1}{2}$, we
conclude
 \begin{equation} \label{eq31}
 \frac{(2^{2\nu(1-\al)} - 2)}
{2^{2\nu(1-\al)} - 1}\big(2^{2\nu(1-\al)}\big)^m \ \leq \ \frac{2}{M} \int_{-M}^M
\left|\frac{K_{\al,\nu}(x+ 2^{-2\nu m -1})-K_{\al,\nu}(x)}{2^{-2\nu
m -1}}\right| dx
 \end{equation}
and \eqref{eq31} equals \eqref{eq10}.

\qed

\section{Singular PDE Solutions of the Aronsson System and the $\infty$-Laplacian.} \label{PDE}

\subsection{Singular Viscosity Solutions of the Aronsson PDE}

We recall here a result established in \cite{K1} by employing the singular function $K_{\al ,\nu}$. We proved that when $n\geq 2$ and $H \in C^1(\R^n)$ is
a Hamiltonian such that some level set contains a line segment, the
Aronsson equation $D^2  u : H_p(Du) \otimes H_p(Du)= 0$ admits explicit entire viscosity solutions. They are superpositions
of a linear part plus a Lipschitz continuous singular part which in general is non-$C^1$ and nowhere twice differentiable. In particular, we supplemented the $C^1$ regularity result of Wang and Yu \cite{W-Y} by deducing that strict level convexity is necessary for $C^1$ regularity of solutions.

\begin{theorem}(cf. \cite{K1})\label{th1} We assume that $H \in
C^1(\R^n)$, $n\geq 2$ and there exists a straight line segment
$[a,b] \sub \R^n$ along which $H$ is constant. Then, for any $F \in W^{1,\infty}_{loc}(\R)$ satisfying $\|F'\|_{L^\infty(\R)}<1$, the formula
 \begin{equation} \label{eq5}
 u(x) \ :=\ \frac{b+a}{2} \cdot x \ +
     \ F\left(\frac{b-a}{2} \cdot x\right),\ \ \
     x \in \R^n,
 \end{equation}
defines an entire viscosity solution $u \in W^{1,\infty}_{loc}(\R^n)$ of the
Aronsson equation \eqref{eq3b}.
\end{theorem}

Here the notation ``$\cdot$" denotes inner product. By employing the particular choice $F:= \int K_{\al ,\nu}$ and variants of this, we provided the following partial answer to the regularity problem:

\begin{corollary}
Strict level convexity of the Hamiltonian $H$ is necessary to obtain
$C^1$ and $C^{1,\be}$ regularity of viscosity solutions to the Aronsson PDE  \eqref{eq3b} in all
dimensions $n\geq 2$.
\end{corollary}

The idea of the proof of Theorem \ref{th1} is the following: suppose first $u$ is smooth. Then, formula \eqref{eq5} is devised in such a way that the image of the gradient $Du$ is contained into the segment $[a,b]$. Hence, $H(Du)$ is constant because $Du(\R^n)$ is contained into a level set of $H$. By rewritting Aronsson's PDE with contracted derivatives, we get
\begin{equation}
A_\infty u \ = \ H_p(Du) \cdot D\big(H(Du)\big)\ = \  0.
\end{equation}
Hence, \eqref{eq5} defines a solution of Aronsson's PDE. However, the previous argument fails when we chose as $F$ the primitive of $K_{\al, \nu}$. For the general case of viscosity solutions, we can use techniques of calculus of the so-called Semijets which are the pointwise generalized derivatives of viscosity solutions to obtain the result. Alternatively, we can use the stability properties of viscosity solutions under limits which claim that local uniform approximation of viscosity solutions produces viscosity solutions to obtain Theorem \ref{th1}.

\subsection{Singular $\infty$-Harmonic Local Diffeomorhisms}

Now we follow \cite{K2} and recall further constructions of singular solutions. 

Let $K\in C^0(\R)$ and define $u : \R^2 \larrow \R^2$ by
\begin{equation}  \label{b8}
u(x,y)\ :=\, \int_0^x e^{i K(t)}dt \ + \ i\int_0^y e^{i K(s)} ds.
\end{equation}
\eqref{b8} defines a $2$-dimensional $\infty$-Harmonic Map, which is singular if $K\not\in C^1(\R)$.

\begin{proposition} (cf. \cite{K2, K6})\label{p2} 
Suppose $\|K\|_{C^0(\R)}< \frac{\pi}{4}$ and let $u$ be given by \eqref{b8}. Then, $u$ is a $C^1(\R^2)^2$-local diffeomorphism and everywhere solution of the PDE system \eqref{inftylap} with contracted derivatives, that is of
\begin{equation}  \label{b9}
Du\, D\left(\frac{1}{2}|Du|^2\right)\ + \ |Du|^2[Du]^\bot \Div\, (Du)\ = \ 0.
\end{equation}
\end{proposition}

In the special case where $u$ is smooth, the main idea is that $u$ equals the sum of two unit speed curves on $\R^2$ in separated variables and as such the Euclidean (Frobenious) matrix norm $|Du|=\sqrt{D_iu_\al D_iu_\al}$ of the gradient is constant. Moreover, in zero-codimension the orthogonal projection $[Du]^\bot$ vanishes identically. Hence, $u$ is a planar $\infty$-Harmonic map. Since the partial derivatives are linearly independent everywhere, the map is both an immersion and a submersion. Hence, by the inverse function theorem it is a local diffeomorphism. 

In the case $K$ is nonsmooth and in particular for $K:=\frac{1}{4}K_{\al, \nu}$, we obtain a $C^{1,\al}(\R^2)^2$ nowhere improvable $\infty$-Harmonic local diffeomorphism which can not be rigorously justified, since the previous argument fails. There appear contracted derivatives which can not be expanded and also multiplications of distributions with nonsmooth functions which are not well-defined. 

As we prove in \cite{K2,K6}, \eqref{b8} can not be rigorously justified as solution classically, strongly, weakly, distributionally or in any other sense. However, $u$ is a Contact Solution of the $\infty$-Laplacian for sufficient exponents ($\al>\frac{1}{2}$). This latter property of our class of functions is essential for the validity of the novel calculus of our Vectorial Contact  Jets (which are the pointwise weak derivatives of the vector case) of the ``weak" theory of Contact Solutions. The latter applies to general fully nonlinear PDE systems.

\subsection{Singular Aronsson Maps}

Again in \cite{K2}, we have proved that under the assumption that some level set of $H$ contains a straight line segment of rank-one matrices and has constant gradient $H_\P$ thereon, there exist singular Aronsson maps, that is, singular solutions of \eqref{aronsyst}. More precisely, given $a,b \in \R^n$, $\eta\in \R^N$ and $f\in C^1(\R)$, we define $u:\R^n \larrow \R^N$ by
\begin{equation}  \label{eq17} 
u(x)\ :=\ \left(x^\top \Big(\frac{b+a}{2}\Big)\right)  \eta \ +\ f\left( x^\top \Big(\frac{b-a}{2}\Big)\right) \eta.
\end{equation}

\begin{proposition} (cf. \cite{K2, K6}) \label{p1}
Let $H\in C^2(\R^N\ot \R^n)$ and suppose that there exist $a,b \in \R^n$, $\eta\in \R^N$, $c\in\R$ and $\textbf{C} \in \R^N\ot \R^n$ such that
\begin{equation}  \label{c1}
[\eta \ot a,\, \eta\ot b] \ \sub \ \{H=c\}\cap\{H_\P =\textbf{C}\}.
\end{equation}
Let $u$ be given by \eqref{eq17} and suppose $\|f'\|_{C^0(\R)}<1$. Then, $u$ is an everywhere solution in $C^1(\R^n)^N$ of the Aronsson PDE system \eqref{aronsyst} contracted, that is of
\begin{equation}  \label{c3}
H_\P(Du)\, D\big(H(Du)\big)\ + \ H(Du)[H_\P(Du)]^\bot \Div \, \big(H_\P(Du)\big)\ = \ 0.
\end{equation}
\end{proposition}

We refrain from presenting more details in this case, since the main ideas relate to the ones of the previous cases.

If $N>1$, \eqref{aronsyst} is a quasilinear nonmonotone system in nondivergence form and generally does not possess classical, strong, weak, measure-valued, distributional or viscosity solutions. In the paper \cite{K6} we introduce our new PDE theory and among other things prove that \eqref{b8}, \eqref{eq17} are appropriately interpreted solutions. 

We emphasize that the specific properties of our class of singular functions imply both the validity of the tools of Viscosity-Contact Solution theories, as well as the inability to rigorously interpret these solutions by means of existing PDE theories in the vector case. In the scalar case, they furnish a regularity counterexample to an important open problem. In the brand new vector case which is at its birth, render the rigorous study of the ``Euler-Lagrange PDEs" of Vector-Valued Calculus of Variations in $L^\infty$ almost impossible without an efficient PDE approach like the one proposed in \cite{K6}.

\medskip
\noindent {\bf Acknowledgement.} {A primitive version of sections 1 and 2 of this work were written when the Author was a doctoral student at the Department of Mathematics, University of Athens, Greece and had been financially supported by the Grant Heraclitus II - \emph{Strengthening human research potential through the implementation of doctoral research}.}

\footnotesize 

\bibliographystyle{amsplain}

\end{document}